         \documentclass[a4paper,12pt,leqno]{article}
\usepackage{amsmath}%
\usepackage{amsfonts}%
\usepackage{amssymb}%
\usepackage{graphicx}
\newtheorem{theorem}{Theorem}

\newenvironment{proof}[1][Proof]{\textbf{#1.} }{\ \rule{0.5em}{0.5em}}
\newtheorem{note}[theorem]{Note}
\begin{document}

\title{Proof of Weierstrass gap theorem
}
\markright{ Weierstrass gap theorem}
\author{Gollakota V V Hemasundar\\ Department of Mathematics, \\SIWS College, Mumbai-400031, INDIA \\E-mail: gvvhemasundar@yahoo.co.in}
\date{}

\maketitle

\begin{abstract}

In this expository note we give  proof of the Weierstrass gap theorem in  Cohomology  terminology .
We analyze gap sequence for finding possible gaps and non-gaps on $X$. 
\end{abstract}
AMS: 30F10

\noindent
\section{Introduction: }
Let $X$ be a compact Riemann surface of genus $g$.
One of the important questions is the existence of meromorphic functions having pole at a single point $p$ on $X$. By  the Riemann Roch theorem, we can show that there always exists a non-constant meromorphic function $f \in \mathcal{M}(X)$ which has a pole of order $\leq g +1$ at $p$ and is holomorphic in $X \setminus \{p\}$. 

One of the basic results in this topic is Weierstrass gap theorem, which states that

\begin{theorem} \label{wgp}\textit{ For a surface of genus $g \geq 1$ there are precisely $g$ integers 
\begin{equation}
1=n_1  <  n_2 < \dots     <n_g <2g
\end{equation}  
such that there does not exist a meromorphic function on $X$ with a pole of order $n_k$ at $p$.} 
\end{theorem}
The numbers $n_j$, for $j = 1, \dots, g$ are called ``gaps'' at $p$ and their complement in $\mathbb{N}$ are called ``non-gaps''.  Further, the sequence is uniquely determined by the point $p$. 

 The case $g=0$ is trivial as there is always a function on the sphere with a simple pole. The case $g=1$ is clear since there is no meromorphic function with a single  simple pole.

In this note we furnish the proof of Theorem \ref{wgp} in the language of sheaf cohomology, in the spirit of proof given in Springer \cite{springer}.  The proof of Noether's theorem was given in Farkas and Kra \cite{farkas} from which the proof of Weierstrass gap theorem was  deduced as a special case.  The proofs are generally given as an application of Riemann-Roch Theorem.  
In section 2, we mention some of the consequences of Riemann-Roch and Serre Duality theorems. The proofs may be found in \cite{forster}.  %
The proof of Weierstrass gap theorem is given in Section 3.   In Section 4, we form a combinatorial problem which appears to be a  by product from the statement of  Weierstrass gap theorem and analyze gaps and non-gaps.

\section{Some Consequences of Riemann-Roch and Serre Duality Theorems}

Let $D$ be a divisor on $X$.

If $deg D < 0$ then $H^0\left(X, \mathcal{O}_D\right) = 0$.

We can use Serre-Duality theorem to obtain equality of dimensions: 
\begin{equation}
\dim H^0\left(X, \Omega_{-D} \right) =  \dim H^1\left(X, \mathcal{O}_D\right)
\end{equation}

For $D=0$ we obtain
\begin{equation}\label{genus of one forms}
\dim H^0\left(X, \Omega \right) = g= \dim H^1\left(X, \mathcal{O}\right)
\end{equation}

Here $H^0\left(X, \Omega \right) = \Omega(X)$ which denotes the sheaf of holomorphic 1-forms on $X$.

The following equality can be obtained as an application of Seree Duality Theorem:

$$
\dim H^1\left(X, \Omega \right) = \dim  H^0\left(X, \mathcal{O} \right)  =1
$$

\begin{note}
This is a known fact again saying there are no non-constant holomorphic maps on a compact Riemann surface.
\end{note}

Suppose $X$ is a compact Riemann surface of genus $g$. Let $K $ be the canonical divisor on $X$. 
Then deg $K = 2 g - 2$.

\section{Proof of Theorem 1 (Weierstrass Gap Theorem)}

\begin{proof}
Suppose $P \in X$.  If $D$ is a zero divisor, then 
$ \dim H^1(X, \mathcal{O}) = g$  and $deg (D)=0$. 

By the Riemann Roch theorem $ \dim H^0(X, \mathcal{O}) = 1$.  Therefore, there are no non-constant holomorphic function on $X$.

 Define the divisor $D_P$ such that 
$$
D_P(P) =
\begin{cases}
1  \mbox{ if } x = P \\
0  \mbox{ if } x \neq P
\end{cases}
$$
$Deg(D_P)=1$.  Once again by the Riemann -Roch theorem, 
$$
\dim H^0(X, \mathcal{O}_{D_P}) = 2-g + \dim H^0(X, \Omega_{-D_P}) 
$$
If $\dim H^0(X, \Omega_{-D_P}) = g$, then $\dim H^0(X, \mathcal{O}_{D_P}) = 2$, hence there exists $f \in \mathcal{M}(X)$ which has a simple pole at $P$ and is holomorphic in $X\setminus \{P\}$. 
 
If $\dim H^0(X, \Omega_{-D_P}) = g-1$, then $\dim H^0(X, \mathcal{O}_{D_P}) = 
1$,  hence there is no meromorphic function which has a simple pole at $P$ and is  holomorphic in $X\setminus \{P\}$.  

Now we want to see the effect of changing $D_P= (n-1)D_P $ to $D_P=nD_P$.  By the Riemann - Roch Theorem

$$
\dim H^0(X, \mathcal{O}_{-D_{(n-1)P}}) = n-g + \dim H^0(X, \Omega_{-D_{(n-1)P}}) 
$$
and 
$$
\dim H^0(X, \mathcal{O}_{-D_{nP}}) = n+1-g + \dim H^0(X, \Omega_{-D_{nP}}) 
$$
If $$\dim H^0(X, \Omega_{-D_{(n-1)P}}) = \dim H^0(X, \Omega_{-D_{nP}}) $$ then

 $$\dim H^0(X, \mathcal{O}_{-D_{nP}}) = \dim H^0(X, \mathcal{O}_{-D_{(n-1)P}}) +1.$$
So there exists a meromorphic function $f \in \mathcal{M}(X)$ with a pole of order $n$ at $P$ and is holomorphic in $X\setminus \{P\}$.

\noindent
If $$\dim H^0(X, \Omega_{-D_{(n-1)P}}) = \dim H^0(X, \Omega_{-D_{nP}})-1 $$ then
$$\dim H^0(X, \mathcal{O}_{-D_{nP}}) = \dim H^0(X, \mathcal{O}_{-D_{(n-1)P}}) .$$ 
So there will not exist a function with a pole of order $n$ at $P$ and holomorphic in  $X\setminus \{P\}$.

So if $\dim H^0(X, \Omega_{-D_{nP}})$ remains the same as $n$ increases by $1$, a new linearly independent function is added in going from  the sheaf $\mathcal{O}_{D_{n_P}}$ to $\mathcal{O}_{D_{(n+1)P}}$.

From the Eq. \ref{genus of one forms} we have $ \dim H^0(X, \Omega) = g $.  We have already seen that if $D$ is the divisor of a non-vanishing meromorphc 1-form on a compact Riemann surface of genus $g$, then deg $(\omega) = 2g-2$.

Let $K$ be its canonical divisor. Then 
 we can deduce that 
 $$\dim H^0(X, \Omega_{-D_{(2g-1)P}})= 0$$. 

 Therefore, the number of times $\dim H^0(X, \Omega_{-D_{nP}})$ does not remain the same must be $g$ times and at each change it decreases by 1. 

It completes the proof.
\end{proof}

\section {Analyzing  gaps and non-gaps}
Suppose $p \in X$. 
If $f$ has a pole of order $s$ at $p$, and $g$ has a pole of order $t$ at $p$, then $fg$ has a pole of order $s+t$ at $p$.  Therefore, the set of non-gaps forms an additive sub-semi group of $\mathbb{N}.$

\noindent
 Let $d$ be the least non-gap value at the point $p$, if $n > d$ is a gap then $n-d$ is again a gap value. Therefore, all  the gaps occur in finite arithmetical sequences  
of the form 
$
 j,  j+ d ,  j + 2d,\dots, j+ \lambda_j d,
$
where $j=1,2, \dots, d-1$ and $\lambda_j= 0,1,2,\dots  .$
See pp. 124 \cite{gunning}.
A point $p$ is called a hyperelliptic Weierstrass point if the non-gap sequence starts with 2 and the hyperelliptic Riemann surfaces are characterized by the gap sequence:
\begin{equation} \label{he}
P= \{1,3, \dots , 2g-1\} 
\end{equation}
hence the non-gaps are
  Q=\{2,4, \dots, 2g\}
	
For the Exceptional Riemann surfaces the gap sequence $P$ and the non-gap sequence $Q$ are given by
\begin{equation} \label{er}
P = \{1,2,3, \dots , g-1, g+1\}  \mbox { and }
Q = \{ g, g+2, \dots , 2g-1\}.
\end{equation}
For example, if   $g=3$ we see the possible gap sequences  are\\
 $\{1,3,5\}$, $\{1, 2,3\}$, $\{1,2,4\}$, $\{1,2,5\}$\\
and corresponding non-gaps are:\\ $\{2,4,6\}$, $\{4, 5,6\}$, $\{3,5,6\}$, $\{3,4,6\}$.


The following problem may be of some interest to see possible gaps and non-gaps at a given point. 



\noindent
\textbf{Question 1:} \\
 Write the numbers  $2$ to $2g-1$ in to two (disjoint) parts $P =\{n_1, n_2, \dots , n_{g-1}\}$ and $Q= \{m_1, m_2, \dots,
 m_{g-1} \}$ such that no  number in $P$ is   a sum of any combination of numbers in $Q$. How many pairs of such P and Q exist?

Clearly then $\{1\} \cup P$ gives possible gaps and $Q \cup \{2g\}$ gives possible non-gaps.

\end{document}